\documentclass[leqno,fleqn]{article}
\usepackage{amsmath}
\usepackage{amssymb}
\usepackage{latexsym}
\usepackage{linguex}
\usepackage{qtree}
\usepackage{graphicx}
\usepackage{theorem}

\usepackage{bussproofs}

\usepackage[backend=biber,bibstyle=ext-authoryear,citestyle=authoryear,articlein=false,dashed=false,giveninits]{biblatex}
\newtheorem{theorem}{Theorem}[section]
\theoremstyle{break}

\newtheorem{lemma}[theorem]{Lemma}
\newtheorem{corollary}[theorem]{Corollary}
\newtheorem{fact}[theorem]{Fact}
{\theorembodyfont{\rmfamily}
\newtheorem{example}[theorem]{Example}
\newtheorem{definition}[theorem]{Definition}
\newtheorem{remark}[theorem]{Remark}}

\newcommand{\be}{\begin{enumerate}}
\newcommand{\bi}{\begin{itemize}}
\newcommand{\ei}{\end{itemize}}
\newcommand{\ee}{\end{enumerate}}
\newcommand{\bqu}{\begin{quote}}
\newcommand{\equ}{\end{quote}}

\newcommand{\lint}{[\![}
\newcommand{\rint}{]\!]}

\newcommand{\proof}{\noindent \emph{Proof. }}
\newcommand{\QED}{\hfill $\Box$ \vspace{3ex}}

\def\emp{\emptyset}

\def\sub{\subseteq}

\def\F{\mathcal{F}}

\def\C{\mathcal{C}}

\def\A{\mathcal{A}}

\def\P{\mathcal{P}}

\def\Q{\mathcal{Q}}

\def\f{\varphi}
\def\p{\psi}

\def\LL{\mathsf{L}}
\def\IL{\mathsf{IPC}}
\def\CL{\mathsf{CPC}}

\title{\Large Carnap's problem for intuitionistic propositional logic}
\author{Haotian Tong\thanks{Tsinghua University (Haotian Tong's work on the paper was done while he was pursuing a second bachelor degree at the Philosophy Department at Tsinghua University, Beijing. His current affiliation is Duke-NUS Medical School, Singapore.
)} \and Dag Westerst\aa hl\thanks{Department of Philosophy at Stockholm University, and Tsinghua-UvA Joint Research Center for Logic, Department of Philosophy, Tsinghua University.}}

\begin{document}

\maketitle

\begin{abstract}
We show that intuitionistic propositional logic is \emph{Carnap categorical}: the only interpretation of the connectives consistent with the intuitionistic consequence relation is the standard interpretation. This holds with respect to the most well-known semantics relative to which intuitionistic logic is sound and complete; among them Kripke semantics, Beth semantics, Dragalin semantics, and topological semantics. It also holds for algebraic semantics, although categoricity in that  case is different in kind from categoricity relative to possible worlds style semantics.
\end{abstract}

\begin{quote}
\end{quote}

\begin{quote}
\end{quote}

\section{Motivation and background\label{s:intro}}

\cite{carnap43} worried about the existence of `non-normal' interpretations of the connectives in classical propositional logic ($\CL$), i.e.\ interpretations that were different from the usual truth tables but still consistent with $\vdash_\CL$. Though neglected for many decades, the issue has been reopened in recent years. One reason was that non-normal interpretations were seen to clash with inferentialist meaning theories (\cite{raatikainen08}). It has been countered that inferentialism is not really threatened by non-normal interpretations (e.g.\ \cite{murzihjortland09}), or that such interpretations are avoided by a proper understanding of inference rules (e.g.\ \cite{murzitopey21}).

In this note we do not discuss the shape or the epistemological status of rules (but see Section \ref{s:discussion} for some discussion), but follow the model-theoretic approach to Carnap's problem initiated in \cite{bonwes16}. That is, a single-conclusion consequence relation $\vdash$ in a logical language is \emph{given} (no matter in what way), and, relative to a semantics which compositionally assigns \emph{semantic values} (set-theoretic objects of various kinds) to formulas, one investigates to what extent the meaning of the logical constants in that language is \emph{fixed} by consistency with $\vdash$. This is a precise model-theoretic task.

A first observation, which, as pointed out in \cite{bonwes16}, is hiding already in Carnap's 1943 book, is that \emph{compositionality} already rules out non-normal interpretations of the connectives in $\CL$, relative to the usual two-valued semantics.\footnote{Compositionality is by now an accepted requirement on formal semantics, but the idea was not around in 1943. The result as stated presupposes that propositional atoms are treated as \emph{variables}, as we do here. Otherwise, one would have to exclude by stipulation the non-standard compositional interpretations making all sentences true, since these would also be consistent with $\vdash_\CL$.
}
More interestingly, classical first-order consequence $\vdash_\mathsf{FOL}$ does \emph{not} by itself fix the meaning of $\forall$, relative to a semantics interpreting quantifier symbols as sets of subsets of the domain (unary generalized quantifiers), but constrains it to be a \emph{principal filter}. But if \emph{permutation invariance} is also required, then the only consistent interpretation is the standard one: $\vdash_\mathsf{FOL}$ forces $\forall x\f(x)$ to mean `for all $x$ in the domain, $\f(x)$ holds'.

Carnap's question can be asked for any consequence relation in any logical language, relative to any formal semantics for that language. \cite{bonwes21} deals with Carnap's problem in modal logic, and \cite{speitel20} discusses it for some logics with generalized quantifiers.

Here we consider Carnap's problem for intuitionistic propositional logic ($\IL$). Relative to which formal semantics should we ask the question? There is by now a plethora of semantics for which $\IL$ is sound and complete: Kripke semantics, Beth semantics, topological semantics, algebraic semantics, \ldots; ideally we would like an answer for each of them.

Actually, there are \emph{two kinds} of semantics for $\IL$: the possible worlds kind and the algebraic kind. The former assigns \emph{subsets of a given domain} as semantic values of formulas. The latter doesn't: semantic values are just elements of the carrier of an algebra. We will see that categoricity means something different in the two cases. In Section \ref{s:categoricity}, we show that \emph{Carnap categoricity} in fact holds for $\IL$, relative to a large number of possible worlds semantics, including the ones just mentioned: only the standard interpretations are consistent with $\vdash_\IL$. In Section \ref{s:algebraicsem} we do the same for algebraic semantics. Section \ref{s:discussion} concludes with discussion and open issues..


Our main reference for the semantics of $\IL$ will be the article \cite{bezhanishvili-holliday19} --- BH19 in what follows --- which surveys and compares a great variety of intuitionistic semantics. In particular, the notion of \emph{nuclear} semantics (introduced in \cite{bezhanishvili-holliday16}) plays a crucial role, as it does in this note. We largely follow the notation and terminology in BH19, and refer to that article for technical notions not explained here.

The propositional language will (with some explicitly noted exceptions) be generated by

\ex.[]
$\f \;:=\; p \;|\; \bot \;|\; \f\land\f \;|\; \f\lor\f \;|\;\f\to\f\;\;\;\;\;(p \in \mathit{Prop})$

where $\mathit{Prop}$ is a denumerable set of propositional variables. So $\neg\f := \f\to\bot$. The single-conclusion consequence relation

\ex.[]
$\vdash_\IL$

in this language is assumed to be familiar, as well as standard intuitionistic Kripke semantics, for which $\vdash_\IL$ is sound and complete.

\section{Possible worlds semantics\label{s:possworldsem}}

What exactly is a possible worlds semantics for $\IL$? Compare the analogous question for basic modal logic. That question has a precise (if initially somewhat surprising) answer: neighborhood semantics. The reason is that, as shown in \cite{bonwes16},  when the semantic values of formulas are subsets of a set $X$ of points (worlds, states, \ldots), classical propositional logic \emph{fixes the meaning} of the propositional connectives: negation is complement, conjunction is intersection, etc.\footnote{In fact, this result is a special case of Theorem \ref{mainthm} below; see Remark \ref{remark:classical}.
}
So the only operator left to interpret is $\Box$, and if no requirements are placed on that interpretation, it must (by compositionality) be a function from sets of points to sets of points. Thus, (local) interpretations can be identified with pairs $(X,F)$, where $F\!:\P(X)\to\P(X)$, and these are exactly the neighborhood frames.\footnote{See \cite{bonwes21}. \cite{pacuit17} is an introduction to neighborhood semantics for modal logic.
}

In the case of $\IL$ there are no further operators to interpret, but surely there should be some constraints on the interpretation of the connectives, since we are after all interested in \emph{intuitionistic} interpretations. Given the great variety of semantics for $\IL$, it is not obvious how to formulate a constraint that fits all.
However, guided by the most well-known instances of such semantics (see below), we suggest, somewhat stipulatively, that an intuitionistic possible worlds semantics assigns \emph{upward closed subsets} (upsets) of partially ordered set (poset) $\F=(X,\leq)$ as semantic values. This immediately rules out the classical interpretation of the connectives (the complement of an upset is usually not an upset).

The elements of the poset can be thought of as \emph{stages}; the rough idea being that information obtained at a stage remains at all later stages. There are numerous implementations and variations of this idea, by Kripke, Dummett, and many others; see the detailed discussion in BH19 and the references therein. But treating semantic values as upsets is common to most of them.

\subsection{Nuclear interpretations}

Let $\mathit{Up}(\F)$ be the set of upsets in $\F=(X,\leq)$, and for $U\sub X$, let $\uparrow\!U$ be the set $\{x\in X\!: \exists y\in U\: y\leq x\}$. Also, for $x\in X$, $\uparrow\!x =\; \uparrow\!\{x\}$; these are the \emph{principal} upsets (the upsets are exactly the unions of principal upsets). 

In algebraic terms, $(\mathit{Up}(\F),\emp,\cap,\cup,\to)$ is a \emph{complete Heyting algebra}, also called a \emph{locale}, where, for $U,V\in \mathit{Up}(\F)$,

\ex.[]
$U\to V = \{x\in X\!:\; \uparrow\!x \cap U \sub V\}$ \footnote{\label{fn:alexandroff}More precisely, it is a \emph{spatial} locale, since it comes from the open sets of a topology; completeness is the property of being closed under arbitrary joins (unions), which holds for all topologies. In this case the topology is rather special (the \emph{Alexandroff} topology): the open sets are the upsets on a poset, and these are also closed under arbitrary meets (intersections). 
}
\vspace{1ex}

\begin{definition}[intuitionistic interpretations and consistency]\label{def:possworldssem}
Let a poset $\F=(X,\leq)$ be given.\vspace{-1ex}
\be
\item An \emph{intuitionistic interpretation of $\IL$ on $\F$} is a map $I^\F$ from $\{\bot,\land,\lor,\to\}$ to functions of the corresponding arity taking upsets to upsets, together with a function $I^{\F,\mathit{at}}\!:\mathit{Up}(\F)\to\mathit{Up}(\F)$ ($\mathit{at}$ for `atom').\vspace{-1ex}
\item A \emph{valuation on $\F$} is a map $v\!:\mathit{Prop}\to \mathit{Up}(\F)$.\footnote{There is no constraint on $v$. But, as we will see, particular intuitionistic interpretations may require all semantic values, even those of propositional atoms, to be a special kind of upset. That's the reason for including the function $I^{\F,\mathit{at}}$ in interpretations; see the next item.
}
\vspace{-1ex}
\item $I^\F$, $I^{\F,\mathit{at}}$, and $v$ compositionally assign a semantic value $\lint\f\rint^{I^\F}_v \in\mathit{Up}(\F)$ of each formula $\f$, as follows:\vspace{-1ex}
\be
\item $\lint p \rint^{I^\F}_v = I^{\F,\mathit{at}}(v(p))$ 
\item $\lint \#(\f_1,\ldots,\f_n) \rint^{I^\F}_v = I^\F(\#)(\lint \f_1 \rint^{I^\F}_v,\ldots,\lint \f_n \rint^{I^\F}_v)$
\ee
\item $I^\F$ is \emph{consistent with a consequence relation $\vdash$} if, whenever $\Gamma\vdash\f$, we have, for all valuations $v$ on $\F$:
\ex.\label{consistencygeneral}
$\bigcap_{\p\in\Gamma}\lint \p \rint^{I^\F}_v \sub\: \lint \f \rint^{I^\F}_v$

We take $\bigcap\emp=X$, so a special case of \ref{consistencygeneral} is that if $\:\vdash\f$, then for all $v$,
\ex.\label{consistencyspecial}
$\lint \f \rint^{I^\F}_v =X$

\ee
\end{definition}

We can now see how familiar semantics for $\IL$ each has a corresponding kind of intuitionistic interpretation, and, among these, a \emph{standard} interpretation. The starting point is the \emph{nuclear} semantic framework from \cite{bezhanishvili-holliday16}.

A \emph{nucleus on $\mathit{Up}(\F)$} is a function $j$ from $\mathit{Up}(\F)$ to $\mathit{Up}(\F)$ such that the following holds for all $U,V \in \mathit{Up}(\F)$:
\be
\item[(i)] $U \sub jU$ \hspace{19ex} (inflationarity)\vspace{-1ex}
\item[(ii)] $jjU \sub jU$ \hspace{17ex} (idempotence) \vspace{-1ex}
\item[(iii)] $j(U\cap V) = jU \cap jV$ \hspace{6ex} (multiplicativity)
\ee
$j$ is \emph{dense} if $j\emp=\emp$. This is an instance of the general notion of a nucleus on a Heyting algebra, wich is defined analogously, with $\cap$ replaced by $\land$, $\emp$ by 0, and $\sub$ by the partial order $a\leq b \Leftrightarrow a\land b= a$. 

It easily follows that $\sub$ in (ii) can be replaced by $=$, and that we have:
\be
\item[(iv)] $U\sub V$ implies $jU \sub jV$ \hspace{2ex} (monotonicity)
\ee
We say that $U\in\mathit{Up}(\F)$ is \emph{fixed}, if it is a fixpoint of $j$, i.e.\ if $jU=U$. Let

\ex.[]
$\mathit{Up}(\F)_j$

be the set of fixed elements of $\mathit{Up}(\F)$. BH19 calls $(\F,j)$ a \emph{nuclear frame}. 

\begin{definition}[nuclear interpretations]\label{def:nuclear}
Let $j$ be a nucleus on $\mathit{Up}(\F)$.\vspace{-1ex}
\be
\item A \emph{nuclear interpretation on} $(\F,j)$ is an intuitionistic interpretation, written $I^{\F,j}$, in which the connectives are interpreted as functions taking fixed upsets to fixed upsets, and $I^{\F,j,\mathit{at}}=j$. Thus, semantic values, calculated as in Definition \ref{def:possworldssem}:3, are fixed upsets. Valuations are as before functions from $\mathit{Prop}$ to $\mathit{Up}(\F)$.
\vspace{-1ex}
\item The \emph{standard} nuclear interpretation, $I^{\F,j}_\text{st}$, is defined as follows, for fixed upsets $U,V$:\vspace{-1ex}
\be
\item $I^{\F,j}_\text{st}(\bot)=j\emp$
\item $I^{\F,j}_\text{st}(\land)(U,V)= U\cap V$
\item $I^{\F,j}_\text{st}(\lor)(U,V)= j(U\cup V)$
\item $I^{\F,j}_\text{st}(\to)(U,V)= \{x\in X\!:\; \uparrow\!x \cap U \sub V\}$
\ee
\ee
\end{definition}

That $I^{\F,j}_\text{st}$ is indeed a nuclear interpretation, i.e.\ that semantic values are fixed upsets, follows from a well-known fact about Heyting algebras: if $j$ is a nucleus on a Heyting algebra $H$, then the algebra $H_j$ of fixpoints of $j$ is again a Heyting algebra (and if $H$ is a locale, so is $H_j$), with $0_j = j0$, $a \land_j b = a\land b$, $a\lor_j b = j(a\lor b)$ (or $\bigvee_jA = j\bigvee \!A$), and $a\to_j b = a\to b$. In particular, if $U,V\in \mathit{Up}(\F)_j$, then $I^{\F,j}_\text{st}(\to)(U,V) \in \mathit{Up}(\F)_j$. Let

\ex.[]
$\mathbf{Up}_j = (\mathit{Up}(\F)_j,\emp_j,\cap_j,\cup_j,\to_j)$\vspace{2ex}

\subsection{Examples}

\begin{example}[Kripke semantics]\label{ex:kripke-int}
Given $\F$, let $j_k$ be the identity function on $\mathit{Up}(\F)$. $j_k$ is obviously a nucleus. So \emph{Kripke frames} are nuclear frames of the form $(\F,j_k)$ (or simply $\F$), and the \emph{standard Kripke interpretation} $I^{\F,j_k}_\text{st}$ is now as in Definition \ref{def:nuclear}, with $j=j_k$. Clearly, these are exactly the usual truth conditions in Kripke semantics for intuitionistic logic. 

More generally, let a \emph{Kripke interpretation on $\F$} be a nuclear interpretation $I^{\F,j_k}$. Semantic values $\lint\f\rint^{I^{\F,j_k}}_v$ are upsets (which are trivially fixed), and $\lint p \rint^{I^{\F,j_k}}_v = j_kv(p) = v(p)$.
\end{example}

\begin{example}[Beth semantics]\label{ex:beth}
There are several versions of Beth semantics, but here we follow BH19: a \emph{Beth frame} is a poset $\F$ (rather than a tree), and a \emph{path} is chain $C$ in $\F$ closed under upper bounds: if $x$ is an upper bound of $C$, then $x\in C$. So the nuclear setting applies without change. Let, for $U\in \mathit{Up}(\F)$, 

\ex.
$j_bU = \{x\in X\!: \text{every path through $x$ intersects $U$}\}$

Then one verfies that $j_b$ is a nucleus on $\mathit{Up}(\F)$ --- the \emph{Beth nucleus} --- and the \emph{standard Beth interpretation} $I^{\F,j_b}_\text{st}$ is as in Definition \ref{def:nuclear} with $j=j_b$. Thus, it differs from the standard Kripke interpretation only in the semantic values of atoms and disjunctions.

In general, we define a \emph{Beth interpretation on $\F$} to be any nuclear interpretation of the form $I^{\F,j_b}$. Semantic values are fixed upsets relative to $j_b$, and $\lint p \rint^{I^{\F,j_b}}_v =  j_bv(p)$.
\end{example}

\begin{example}[Dragalin semantics]
Dragalin generalized Beth semantics by considering a more general notion of a path, called a \emph{development} by Bezhanishvili and Holliday. Each $x\in X$ is a associated with a set $D(x)$ of subsets of $X$ satisfying certain conditions. Dragalin proved that

\ex.[]
$j_DU = \{x\in X\!: \text{every development in $D(x)$ intersects $U$}\}$

is a nucleus on $\mathit{Up}(\F)$, so Dragalin semantics is an instance of nuclear semantics. The importance of this kind of nuclear semantics is seen from the result in \cite{bezhanishvili-holliday16} that \emph{every} locale can be realized as the set of fixed upsets of a nuclear frame of the form $(\F,j_D)$ (this does not hold for Beth semantics).

Just as before, one has the \emph{standard Dragalin interpretation}, and a notion of an arbitrary \emph{Dragalin interpretation}.
\end{example}

We have now seen three main kinds of what we may call \emph{nuclear semantics}, in a sense that can be made precise as follows:

\begin{definition}[nuclear semantics]\label{def:nuclearsemantics}
A \emph{nuclear semantics} for $\IL$ is a class $\C$ of nuclear frames. For $(\F,j)\in\C$, a \emph{$\C$-interpretation on $(\F,j)$} is a nuclear interpretation $I^{\F,j}$.
\end{definition}

For example, \emph{Beth semantics} is then identified with the class of nuclear frames $(\F,j_b)$ where $\F$ is a Beth frame, and \emph{Beth interpretations on} $(\F,j_b)$ are as in Example \ref{ex:beth} above. BH19 shows that $\vdash_\IL$ is sound for \emph{any} nuclear semantics, if we use the standard interpretation. The next semantic framework for $\IL$ is not strictly nuclear, but almost.

\begin{example}[topological semantics]\label{ex:topo-int}
Topological semantics is the oldest formal semantics for $\IL$, going back to \cite{stone37} and \cite{tarski38}. That it can be construed as an instance of nuclear semantics follows from a special case of a theorem of Dragalin; for a proof of the general case in the present setting see \cite{bezhanishvili-holliday16} (Theorem 2.8).

Let $\Omega(X)$ be a topology on $X$, i.e.\ $\Omega(X)$ is its set of opens.\footnote{So $\Omega(X)\sub\P(X)$, $\emp,X\in \Omega(X)$, and $\Omega(X)$ is closed under finite intersections and arbitrary unions.
}
Then $\mathbf{\Omega}(X) = (\Omega(X),\emp,\cap,\cup,\to)$ is a complete Heyting algebra, with

\ex.[]
$U\to V = \mathit{int}((X\!-\!U)\cup V)$

where $\mathit{int}$ is the \emph{interior operation}: for $Y\sub X$,

\ex.[]
$\mathit{int}(Y) = \bigcup\{U\!\in\Omega(X)\!: \:U\sub Y\}$

The open sets themselves need not be upsets of a partial order, but Dragalin's result is that $\mathbf{\Omega}(X)$ is  \emph{isomorphic} to the fixed upsets of a nuclear frame.

Let $\F$ be the poset $(\Omega(X)^-,\supseteq)$, where $\Omega(X)^- = \Omega(X) - \{\emp\}$ (so the upsets of $\F$ are the downsets of $(\Omega(X)^-,\sub)$). We use $A,B,\ldots$ for upsets of $\F$, and $Z,U,V\ldots$ for open sets. Define:

\ex.[]
$jA = \;\uparrow\! \bigcup A$

\begin{theorem}[Dragalin]\label{dragalinthm}
$j$ is a dense nucleus on $\F$, and the function $h(U) =\; \uparrow\!U$ is an isomorphism from $\mathbf{\Omega}(X)$ to $\mathbf{Up}(\F)_j$.
\end{theorem}

So in this sense, a topology can be seen as a nuclear frame, defined via the bijection $h$. Next, the \emph{standard topological semantics} on $\Omega(X)$ for $\IL$, originating from Stone and Tarski, is given by the interpretation we write $I^{\Omega(X)}_\text{st}$, as follows. Valuations are functions from \emph{Prop} to $\Omega(X)$, and semantic values are calculated with the these functions:

\ex.\label{topotruth}
\a.\label{topotrutha} $I^{\Omega(X)}_\text{st}(\bot)=\emp$\vspace{1ex}
\b. $I^{\Omega(X)}_\text{st}(\land)(U,V)= U\cap V$\vspace{1ex}
\b. $I^{\Omega(X)}_\text{st}(\lor)(U,V)= U\cup V$\vspace{1ex}
\b. $I^{\Omega(X)}_\text{st}(\to)(U,V)= \mathit{int}((X\!-\!U)\cup V)$

This is not an intuitionistic possible worlds semantics in the sense of Definition \ref{def:possworldssem}, but it is isomorphic to one:

\begin{lemma}\label{topohelplemma}
If $\F$, $j$, and $h$ are as in Theorem \ref{dragalinthm}, then, for all $A,B\in\mathit{Up}(\F)_j$, $I^{\F,j}_\text{st}(\land)(A,B) = h(I^{\Omega(X)}_\text{st}\!(\land)(h^{-1}(A),h^{-1}(B)))$, and similarly for the other connectives.
\end{lemma}

\proof
The proof is a routine verification, using the density of $j$ for $\bot$, and, for $\lor$ and $\to$, the following observation from the proof of Theorem \ref{dragalinthm}:

\ex.\label{jA=Afact}
$jA=A\;$ iff $\;A=\emp$ or $A=h(U)$ for some $U\in \Omega(X)^-$.

We omit the details.
\QED







In analogy with the preceding examples, let us say that a \emph{topological interpretation of $\IL$ on} $\Omega(X)$ is a map $I^{\Omega(X)}$ from $\{\bot,\land,\lor,\to\}$ to functions of the corresponding arity taking open sets to open sets. A \emph{topological valuation} is a map $v\!: \mathit{Prop} \to \Omega(X)$, and semantic values $\lint\f\rint^{I^{\Omega(X)}}_v$ are assigned compositionally just as before. Also, consistency of such an interpretation with a consequence relation is defined as before.
\end{example}

\section{Carnap categoricity\label{s:categoricity}}

$\vdash_\IL$ is sound and complete for all the semantics (classes of nuclear frames, or topological spaces) discussed in the preceding section --- relative, of course, to what we are here calling the standard interpretations; see BH19 for proofs. However, our interest is not completeness, but categoricity. For each particular semantics $\C$ mentioned above, we defined the notion of a (local) $\C$-interpretation of the connectives, and described how semantic values relative to such an interpretation are computed for each formula. 

The standard interpretations are all consistent with $\vdash_\IL$. But could there be \emph{other} interpretations, of the relevant kind, that are \emph{also} consistent with $\vdash_\IL$? This is Carnap's question. 

\begin{definition}[categoricity]\label{def:categoricity}
A consequence relation $\vdash$ in the propositional language is (Carnap) \emph{categorical} with respect to a nuclear semantics $\C$ if, for every $(\F,j)\in\C$, the only $\C$-interpretation on $(\F,j)$ consistent with $\vdash$ is $I^{\F,j}_\text{st}$. Similarly for topological semantics. If $\vdash$ is associated with a particular logic $\LL$, we also say that $\LL$ is \emph{categorical} (with respect to $\C$) when $\vdash$ is.
\end{definition}

Our main result is the following.

\begin{theorem}\label{mainthm}
Let $(\F,j)$ be a nuclear frame, and $I^{\F,j}$ a nuclear interpretation which is consistent with $\vdash_\IL$. Then $I^{\F,j} = I^{\F,j}_\text{st}$.
\end{theorem}

\proof
We first show:

\ex.[(a)]
$I^{\F,j}(\land)$ is standard.

Here and below, $U,V\in \mathit{Up}(\F)_j$. Since $p\land q\vdash_\IL p$, we have $I^{\F,j}(\land)(U,V)\sub U$. Let us see in detail how this happens. Let $v$ be a valuation such that $v(p)=U$ and $v(q)=V$. The requirement that $v(p),v(q)\in \mathit{Up}(\F)$ is satisfied. By the truth definition (Definition \ref{def:possworldssem}:3) and consistency with $\vdash_\IL$, we have 
\[
\lint p\land q\rint^{I^{\F,j}}_v \!= I^{\F,j}\!(\land)(\lint p\rint^{I^{\F,j}}_v\!,\lint q\rint^{I^{\F,j}}_v) = I^{\F,j}\!(\land)(jU,jV) = 
I^{\F,j}(\land)(U,V) \sub \lint p\rint^{I^{\F,j}}_v \!= U
\]
(note that $U$ and $V$ are fixpoints). Similarly, $I^{\F,j}(\land)(U,V)\sub V$. In the other direction, we use the fact that $p,q\vdash_\IL p\land q$ to see that $U\cap V \sub I^{\F,j}(\land)(U,V)$. Thus, $I^{\F,j}(\land)(U,V) = U\cap V = I^{\F,j}_\text{st}(\land)(U,V)$. Next,

\ex.[(b)]
$I^{\F,j}(\bot)$ is standard.

This follows from the fact that $\bot\vdash_\IL p$: letting $v(p) = \emp$, we obtain that $I^{\F,j}(\bot) \sub \lint p\rint^{\F,j}_v = j\emp$. Since $\emp \sub I^{\F,j}(\bot)$ we also have $j\emp \sub jI^{\F,j}(\bot) = I^{\F,j}(\bot)$, by the monotonicity of $j$ and the assumption that $I^{\F,j}(\bot)$ is a fixed upset. Thus, $I^{\F,j}(\bot) = j\emp = I^{\F,j}_\text{st}(\bot)$. 

We next establish some facts about the interpretation of $\to$.

\ex.[(c)]
$I^{\F,j}(\to)(U,V) \:\sub\: I^{\F,j}_\text{st}(\to)(U,V)$

To see this, take $x\in I^{\F,j}(\to)(U,V)$. In order to show that $x\in I^{\F,j}_\text{st}(\to)(U,V) = I^{\F}_\text{st}(\to)(U,V)$, we must show $\uparrow\!x \cap\: U \sub V$. From the fact that $p,p\to q\vdash_\IL q$, we obtain, as should now be clear,

\ex.[]
$U \cap I^{\F,j}(\to)(U,V) \:\sub\: V$

Take $y\in\; \uparrow\!\!x \cap\: U$. Since $y\geq x$ and $I^{\F,j}(\to)(U,V)$ is an upset, $y\in I^{\F,j}(\to)(U,V)$. So $y\in V$ follows, and (c) is proved.

\ex.[(d)]
$U\sub V\;$ iff $\;I^{\F,j}(\to)(U,V) = X$

Suppose $U\sub V$, so $U\cap V = U$. Since $\vdash_\IL p\land q \to q$, we obtain, from consistency with $\vdash_\IL$ and the fact that $\land$ is standard, that $I^{\F,j}(\to)(U\cap V,V)=X$, i.e.\ that $I^{\F,j}(\to)(U,V)=X$. In the other direction, suppose $I^{\F,j}(\to)(U,V)=X$. By (c), $I^{\F,j}_\text{st}(\to)(U,V)=X$. Then for all $x\in X$, $x\in U \Rightarrow x\in V$ (since $x\leq x$), i.e.\ $U\sub V$. This proves (d). 

We are now able to show:

\ex.[(e)]
$I^{\F,j}(\lor)$ is standard.

From $p\vdash_\IL p\lor q$ and $q\vdash_\IL p\lor q$ we obtain that $U\cup V\sub I^{\F,j}(\lor)(U,V)$. By monotonicity and since $I^{\F,j}(\lor)(U,V)$ is a fixed upset, we have $j(U\cup V) \sub jI^{\F,j}(\lor)(U,V) = I^{\F,j}(\lor)(U,V)$. In the other direction, let $v$ be a valuation such that $v(p)=U$, $v(q)=V$, and $v(r)=U\cup V$. (Note that $U\cup V$ is an upset, even though it need not be fixed.\footnote{But we could equally well have taken $v(r)=j(U\cup V)$.
})
Since $U \sub U\cup V$, we have $jU \sub j(U\cup V)$, and so by (d), $I^{\F,j}(\to)(jU,j(U\cup V)) = I^{\F,j}(\to)(U,j(U\cup V)) = \lint p\to r\rint^{I^{\F,j}}_v = X$. Similarly, $\lint q\to r\rint^{I^{\F,j}}_v = X$. Then, from consistency and the fact that

\ex.[]
$p\to r,q\to r,p\lor q\vdash_\IL r$

we obtain $I^{\F,j}(\lor)(U,V) \sub j(U\cup V)$. That is, $I^{\F,j}(\lor)(U,V) = j(U\cup V) = I^{\F,j}_\text{st}(\lor)(U,V)$, and (e) is proved.

Finally, we show:

\ex.[(f)]
$I^{\F,j}(\to)$ is standard.

By (c), it suffices to show that $I^{\F,j}_\text{st}(\to)(U,V) \sub I^{\F,j}(\to)(U,V)$. Thus, take $x\in I^{\F,j}_\text{st}(\to)(U,V)$, i.e.\ $\uparrow\!x\cap U \sub V$. Using the multiplicativity and monotonicity of $j$, we have

\ex.[]
$j(\uparrow\!x\cap U) = j \!\!\uparrow\!x \cap jU = j \!\!\uparrow\!x \cap U \sub jV = V$

Let $v$ be such that $v(p) = \;\uparrow\!\!x$, $v(q)=U$, and $v(r)=V$. Thus, since $\land$ is standard,

\ex.[]
$\lint p\land q\rint^{I^{\F,j}}_v =\: j \!\!\uparrow\!x \cap U \sub V =\, \lint r\rint^{I^{\F,j}}_v$

By (d), it follows that $\lint p\land q\to r\rint^{I^{\F,j}}_v =\, X$. Since

\ex.[]
$p,p\land q\to r \vdash_\IL q\to r$

we conclude that $ j \!\!\uparrow\!\!x \sub I^{\F,j}(\to)(U,V)$. And since $x\in \:\uparrow\!\!x \sub j\!\!\uparrow\!\!x$ (by inflationarity), we have that $x\in I^{\F,j}(\to)(U,V)$. That is, $I^{\F,j}(\to)(U,V) = I^{\F,j}_\text{st}(\to)(U,V)$. 

This concludes the proof that $I^{\F,j} = \,I^{\F,j}_\text{st}$.
\QED

\begin{corollary}
$\IL$ is Carnap categorical with respect to any nuclear semantics, such as Kripke semantics, Beth semantics, or Dragalin semantics. It is also categorical with respect to topological semantics. Thus, the only interpretation, of the respective kind, of the propositional connectives that is consistent with $\vdash_\IL$ is the standard interpretation.
\end{corollary}

\proof
This is immediate from the theorem for nuclear semantics. For topological semantics, let $\Omega(X)$ be a topology on $X$, and suppose the topological interpretation $I^{\Omega(X)}$ is consistent with $\vdash_\IL$. Let $\F$, $j$, and $h$ be as in Theorem \ref{dragalinthm}. Define a nuclear interpretation $I^{\F,j}$ as follows: $I^{\F,j}(\land)(A,B) = h(I^{\Omega(X)}(\land)(h^{-1}(A),h^{-1}(B)))$, and similarly for the other connectives. Also, if $v\!:\mathit{Prop}\to \mathit{Up}(\F)$, define the topological valuation $v^*\!:\mathit{Prop}\to \Omega(X)$ by $v^*(p)=h^{-1}(v(p))$. It follows by an easy induction over formulas (using \ref{jA=Afact} for the atomic case) that for every $\f$,

\ex.[]
$\lint\f\rint^{I^{\F,j}}_{v} = \;h(\lint\f\rint^{I^{\Omega(X)}}_{v^*})$

Since $h$ is monotone ($U\sub V \Rightarrow h(U)\sub h(V)$), this entails that $I^{\F,j}$ is consistent with $\vdash_\IL$. Then we have: $I^{\Omega(X)}(\land)(U,V) = h^{-1}(I^{\F,j}(\land)(h(U),h(V)))$ (by definition) $= h^{-1}(I^{\F,j}_\text{st}(\land)(h(U),h(V)))$ (by the theorem) $= I^{\Omega(X)}_\text{st}(\land)(U,V)$ (by Lemma \ref{topohelplemma}). Similarly for the other connectives. Thus, $I^{\Omega(X)} =\, I^{\Omega(X)}_\text{st}$.
\QED

\begin{remark}\label{remark:classical}
The Carnap categoricity of $\CL$ relative to classical possible worlds semantics, proved in \cite{bonwes16}, is a special case of the result for $\IL$, since on posets that are sets of isolated points, $\IL$ and $\CL$ coincide. In more detail: suppose $W$ is any non-empty set (of `worlds') and $I^W$ an interpretation of the connectives over $W$ --- i.e.\ $I^W$ assigns to each connective a function on $\P(W)$ of appropriate arity --- which is consistent with $\vdash_\CL$. The upsets of the poset $\F = (W,\{(x,x)\!: x\in W\})$ are exactly the subsets of $W$. Thus $I^\F = I^W$ is in fact a Kripke interpretation, in the sense of Example \ref{ex:kripke-int}, which is consistent with $\vdash_\CL$, hence with $\vdash_\IL$. By the theorem, $I^\F$ is standard, and so $I^W$ is standard (for example, $I^W(\to)(U,V)=(W\!-\!U)\cup V$).

Similarly, letting $\F$ be a single reflexive point $(\{x\},\{(x,x)\})$, the categoricity of $\CL$ relative to classical 2-valued semantics, mentioned in Section \ref{s:intro}, also follows from Theorem \ref{mainthm}: in this case the interpretation functions are truth functions, and the truth values 0 and 1 correspond to the upsets $\emp$ and $\{x\}$.
\end{remark}

As the remark illustrates, the power of the intuitionistic categoricity theorem comes from its strictly local character: the result holds for \emph{each} nuclear frame, no matter how trivial.

The proof of Theorem \ref{mainthm} uses essentially the assumption that consistency with $\vdash_\IL$ means that \ref{consistencygeneral} holds (for all $v$). One might be tempted to suppose that since $\IL$ satisfies the Deduction Theorem, it would suffice that all $\IL$-\emph{theorems} are valid (so that \ref{consistencyspecial} holds for all $v$). It is instructive to see why this is \emph{not} the case. 

To begin, we can see where the proof breaks down. Suppose $\p\vdash_\IL\f$. By the Deduction Theorem, we have $\vdash_\IL\f\to\p$, so by the soundness assumption,  $\lint\p\to\f\rint^{I^{\F,j}}_v=I^{\F,j}(\to)(\lint\p\rint^{I^{\F,j}}_v,\lint\f\rint^{I^{\F,j}}_v)=X$ for all $v$. We wish to conclude that $\lint\p\rint^{I^{\F,j}}_v \sub \lint\f\rint^{I^{\F,j}}_v$. This would follow from (d) in the proof, but (d) presupposes consistency in the stronger sense. More precisely, (d) relies on (c), whose proof in turn uses that from $p,p\to q\vdash_\IL q$ we are able to conclude that, for any fixed upsets $U,V$, we have $U \cap I^{\F,j}(\to)(U,V) \sub V$. But it should be fairly clear that from $\vdash_\IL p\to((p\to q)\to q)$, no such conclusion can be drawn.

In fact, Wesley Holliday found a counter-example (\emph{p.c}): a non-standard interpretation of the connectives which validates all $\IL$-theorems, and hence is not consistent with $\vdash_\IL$ in the sense of Definition \ref{def:possworldssem}:4. With his kind permission, we present the example here.\footnote{The example was originally intended to show that the stronger notion of consistency is needed for the result in \cite{bonwes16} about $\CL$; see Remark \ref{remark:classical}. Here we have adapted it to $\IL$.
}

\begin{example}[Holliday]\label{ex:holliday}
Let $\Omega(X)$ be any topological space. We will define a topological interpretation $I^{\Omega(X)}$ (see Example \ref{ex:topo-int}) such that (a) if $\vdash_\IL\f$ then for all valuations $v$ on $\Omega(X)$, $\lint\f\rint^{I^{\Omega(X)}}_v = X$, but (b) $I^{\Omega(X)}\neq I^{\Omega(X)}_\text{st}$. This gives a generic topological counter-example. For more concreteness, we can start with a Kripke frame $\F=(X,\leq)$ and let $\Omega(X)$ be the Alexandroff topology described in footnote \ref{fn:alexandroff}. As made clear in BH19, Sect.\ 2.3, Kripke semantics essentially \emph{is} topological semantics based on Alexandroff spaces. It is easy to see that the counter-example then becomes a non-standard Kripke interpretation $I^{\F,j_k}$ (see Example \ref{ex:kripke-int}) which validates all theorems of $\IL$.

For this example it is easier to use $\neg,\land,\lor,\to$ as primitive connectives. \ref{topotrutha} is then replaced by

\ex.[\ref{topotrutha}$'$]
$I^{\Omega(X)}_\text{st}(\neg)(U)= \mathit{int}(X-U)$

Recall the \emph{closure} operation $\mathit{cl}$, the dual of $\mathit{int}$: $\mathit{cl}(Y)$ is the smallest closed set (set whose complement is open) containing $Y$. An easy calculation shows

\ex.\label{Wes1}
$\lint\neg\neg\f\rint^{I^{\Omega(X)}_\text{st}}_v = \mathit{int}(\mathit{cl}(\lint\f\rint^{I^{\Omega(X)}_\text{st}}_v))$

and thus

\ex.\label{Wes2}
$\lint\neg\neg\neg\f\rint^{I^{\Omega(X)}_\text{st}}_v = \mathit{int}(X-\mathit{int}(\mathit{cl}(\lint\f\rint^{I^{\Omega(X)}_\text{st}}_v)))$

Now define $I^{\Omega(X)}$ as follows.

\ex.\label{Wes3}
\a. $I^{\Omega(X)}(\land)(U,V) = \mathit{int}(\mathit{cl}(U)) \cap \mathit{int}(\mathit{cl}(V))$\vspace{1ex}
\b. $I^{\Omega(X)}(\neg)(U) = \mathit{int}(X-\mathit{int}(\mathit{cl}(U)))$\vspace{1ex}
\b. $I^{\Omega(X)}(\to)(U,V) = I^{\Omega(X)}(\neg)(I^{\Omega(X)}(\land)(I^{\Omega(X)}(\neg)(U),V))$\vspace{1ex}
\b. $I^{\Omega(X)}(\lor)(U,V) = I^{\Omega(X)}(\neg)(I^{\Omega(X)}(\land)(I^{\Omega(X)}(\neg)(U),I^{\Omega(X)}(\neg)(V)))$

Thus, we are putting double negations in front of negated formulas and the conjuncts of conjunctions, whereas $\to$ and $\lor$ are defined classically from $\land$ and $\neg$. Clearly, $I^{\Omega(X)}$ is non-standard: for example, just find a space with open sets $U,V$ such that $\mathit{int}(\mathit{cl}(U)) \cap \mathit{int}(\mathit{cl}(V)) \neq U\cap V$.

Now we use a particular \emph{negative translation} of classical into intuitionistic propositional logic, defined as follows:

\ex.\label{Wes4}
\a. $g(p)=p$
\b. $g(\f\land\p) = \neg\neg g(\f) \land \neg\neg g(\p)$
\b. $g(\neg\f) = \neg\neg\neg g(\p)$
\b. $g(\f\to\p) = \neg(g(\f) \land \neg g(\p))$
\b. $g(\f\lor\p) = \neg(\neg g(\f) \land \neg g(\p))$

Using well-known facts about negative translations, it is not hard to show that

\ex.\label{Wes5}
$\models_\CL\f\;$ iff $\;\vdash_\IL g(\f)$.\footnote{$g$ is similar to the \emph{Kolmogorov translation}, say, $G$, which puts a double negation in front of all subformulas. One can show that if $\f$ is not an atom, $\vdash_\IL \!g(\f) \leftrightarrow G(\f)$. It is well-known that $\vdash_\CL\f \Leftrightarrow \;\,\vdash_\IL G(\f)$, and since no atom is a theorem, \ref{Wes5} follows. \cite{ferreira12} is a survey of various negative translations.
}

Finally, observe that $I^{\Omega(X)}$ interprets a formula $\f$ just as $g(\f)$ is standardly interpreted in topological semantics. In other words, for each topological valuation $v$ we have:

\ex.\label{Wes6}
$\lint\f\rint^{I^{\Omega(X)}}_v\! = \lint g(\f)\rint^{I^{\Omega(X)}_\text{st}}_v$

This is proved by a straightforward inductive argument, using the standard topological truth definition and \ref{Wes1} -- \ref{Wes4}. Thus, if $\vdash_\IL\f$, then $\vdash_\CL\f$, hence $\vdash_\IL g(\f)$ by \ref{Wes5}, and so for any topological valuation $v$, $\lint g(\f)\rint^{I^{\Omega(X)}_\text{st}}_v = X$, which by \ref{Wes6} entails that $\lint\f\rint^{I^{\Omega(X)}}_v=X$. That is, $I^{\Omega(X)}$ validates all $\IL$ theorems.\footnote{Indeed, it validates all $\CL$ theorems, by \ref{Wes5}, and since the above argument in fact works for all valuations $v\!: \mathit{Prop} \to \P(X)$. As Wesley Holliday pointed out, $I^{\Omega(X)}$ evaluates complex formulas using the \emph{double negation nucleus}, i.e.\ in the Heyting algebra of \emph{regular open} sets (sets $U$ such that $U = \mathit{int}(\mathit{cl}(U)))$, which is a Boolean algebra.
}

On the other hand, it is easy to see that $I^{\Omega(X)}$ is not consistent with $\vdash_\IL$. For example, one can have $\lint p\rint^{I^{\Omega(X)}}_v \!= \lint p\to q\rint^{I^{\Omega(X)}}_v \!= X$, while $\lint q\rint^{I^{\Omega(X)}}_v\!\neq X$. 
\end{example}

\begin{remark}\label{remark:classicalML}
It may be worth noting that, by contrast, categoricity facts about classical \emph{modal} logic (see \cite{bonwes21}), only require validating all theorems. This is precisely because it is classical. If $\mathsf{L}$ is a modal logic, $\Gamma\vdash_\mathsf{L}\f$ means by definition that for some $\p_1,\ldots,\p_n\in\Gamma$, we have $\vdash_\mathsf{L}\p_1\land\ldots\land\p_n\to\f$. Thus, if it holds for all $v$ that $\lint\p_1\land\ldots\land\p_n\to\f\rint^{(X,F)}_v=X$, we can conclude, by the fact that $\land$ and $\to$ are standard (which, as we saw in Remark \ref{remark:classical}, follows from Theorem \ref{mainthm}), that $\bigcap_{i=1}^n\lint\p_i\rint^{(X,F)}_v\sub \lint\f\rint^{(X,F)}_v$, and hence that $\bigcap_{\p\in\Gamma}\lint\p\rint^{(X,F)}_v\sub \lint\f\rint^{(X,F)}_v$.
\end{remark}

\section{Algebraic semantics\label{s:algebraicsem}}

Algebraic semantics for $\IL$ standardly interprets the connectives in a Heyting algebra $\A = (A,0^\A,1^\A,\land^\A,\lor^\A,\to^\A)$: a valuation $v$ maps atoms to $A$, and semantic values $\lint\f\rint^\A_v$ in $A$ are given by the unique homomorphism $\overline{v}$ from the syntax algebra to $\A$ that extends $v$.\footnote{In this section we take 1 ($\top$) to be in the signature (as is commonly done).
}
Bezhanishvili and Holliday discuss the view that algebraic semantics is just `syntax in disguise' (quoted from \cite{vanbenthem01}, p.\ 358), and that algebraic completeness of $\IL$ via the usual Tarski-Lindenbaum construction isn't really illuminating about the intuitionistic meaning of the connectives (BH19, end of sect.\ 2.1). They counter that when Heyting algebras are defined order-theoretically, and $\leq$ is seen as entailment, completeness for Heyting algebras is in fact quite illuminating.\footnote{One might also note that completeness for `concrete' semantics, like Kripke semantics or topological semantics, \emph{follows} from algebraic completeness via well-known representation theorems. Arguably, however, the real work here lies in the representations theorems, not in the Lindenbaum algebras of logical systems.
}

From the categoricity perspective, Carnap's question for algebraic semantics is conceptually somewhat different from the case of possible worlds semantics. In the latter case, given, say, a nuclear frame, there are many different putative interpretations of the connectives \emph{on that frame}, and we ask which ones are consistent with $\vdash_\IL$. Relative to an algebra $\A$, on the other hand, the interpretation of $\to$, for example, simply \emph{is} $\to^\A$. That is, \emph{the algebra itself} is the interpretation. Indeed, we shall see that there are sublogics of $\IL$ which are categorical with respect to algebraic semantics, but not with respect to Kripke semantics.

We ask, then: Are there \emph{other} algebras than Heyting algebras --- the agreed-on standard interpretations --- consistent with $\vdash_\IL$? It seems plausible that the answer is No: consistency with $\vdash_\IL$ forces the algebra to be Heyting.

To make the question meaningful we should specify a class of algebras among which the Heyting algebras can be singled out. It has to be algebras for which a suitable notion of consistency with a consequence relation makes sense. For the record, we now describe how this can be done. We hasten to add, however, that when it comes to categoricity rather than completeness, the claim that algebraic semantics is `syntax in disguise' seems rather convincing. This will be clear from the proof of Theorem \ref{mainthmalg} below. But since Carnap's question in this algebraic context is new (as far as we know), we shall spell out the fairly obvious answer.

\subsection{Algebraic interpretations and consistency}

The syntax algebra of propositional logic is a \emph{term algebra} with countably many generators, which means that for \emph{any} algebra $\A$ of the same signature, every map $v$ from propositional atoms to $A$ extends to a unique homomorphism $\overline{v}$ to $\A$. But it would make no sense to let every algebra of that signature be a putative interpretation of the connectives. Interpretations are relative to a \emph{semantics}, and a semantics needs a notion of \emph{truth}, or of \emph{entailment}. We could of course stipulate that $\overline{v}(\f)=1^\A$ means that $\f$ is true in $\A$ under $v$. But nothing can be done with that stipulation unless we know more about the role of $1^\A$. For example, why not $0^\A$ instead?

On the other hand, we want to make as few assumptions as possible. Without further ado, here is a suggestion.

\begin{definition}\label{def:algsem}
$\A = (A,0^\A,1^\A,\land^\A,\lor^\A,\to^\A)$ is an \emph{algebraic interpretation} if $\leq^\A$ defined by $a\leq^\A b \Leftrightarrow a\land^\A b=a$ is a partial order with $1^\A$ as its largest element.
\end{definition}

Then we define consistency with a consequence relation $\vdash$ as follows.\footnote{We are now in the framework that \cite{humberstone11} calls \emph{$\leq$-based algebraic semantics}. The definition is suggested by his Remark 2.14.7(i).
}

\begin{definition}\label{def:algcons}
An algebraic interpretation $\A$ is \emph{consistent with} $\vdash$, if, whenever $\p_1,\ldots,\p_n\vdash\f$ holds, we have for all valuations $v$ on $A$ and all $c\in A$, that if $c \leq^\A \overline{v}(\p_i)$ for $i=1,\ldots,n$, then $c \leq^\A \overline{v}(\f)$. 
\end{definition}

If there are no assumptions, i.e.\ if $\vdash\!\f$, the antecedent of the requirement is vacuously satisfied, and the consequent becomes $\overline{v}(\f)=1^\A$. It may seem more natural to require that if $c$ is the greatest lower bound of $\overline{v}(\p_1),\ldots,\overline{v}(\p_n)$, then $c \leq^\A \overline{v}(\f)$, and indeed this is equivalent, provided glb's \emph{exist}. In that case, $\mathit{glb}^\A\{a_1,\ldots,a_n\}=a_1\land^\A\ldots\land^\A a_n$, and the consistency requirement would simply be that $\p\vdash\f$ entails $\overline{v}(\p) \leq^\A \overline{v}(\f)$ for all valuations $v$ on $A$.\footnote{\label{fn:consistency}The consistency requirement only concerns the reduct $(A,1^\A,\land^\A)$; $1^\A$ is not explicitly mentioned. Generalizing, for algebras without a top element, we can replace ``$\overline{v}(\f)=1^\A$'' with ``for all $\theta$, $\overline{v}(\theta) \leq^\A \overline{v}(\f)$'' in order to deal with theorems. Logics whose set of theorems coincides with the set of formulas derivable from every formula are sometimes called \emph{non-pseudo-axiomatic}.
}
But Definition \ref{def:algsem} doesn't require glb's to exist, and the next result shows that such a requirement is not necessary.

\subsection{Algebraic categoricity}

We now have a precise formulation of Carnap's question for algebraic semantics. The next fact gives the answer, here formulated for $\CL$ as well. Let $\mathit{HA}$ ($\mathit{BA}$) be the class of Heyting (Boolean) algebras.

\begin{theorem}\label{mainthmalg}
Let $\A$ be an algebraic interpretation.\vspace{-1ex}
\be
\item $\A$ is consistent with $\vdash_\IL$ iff $\A\in\mathit{HA}$.\vspace{-1ex}
\item $\A$ is consistent with $\vdash_\CL$ iff $(A,0^\A,1^\A,\land^\A,\lor^\A,'^{\A})\in\mathit{BA}$, where $a'^{\A} := a\to^\A0^\A$.
\ee
\end{theorem}

\proof
The right-to-left directions are just the soundness of $\IL$ ($\CL$) for the class of Heyting (Boolean) algebras. In the other direction, suppose $\A$ is an algebraic interpretation consistent with $\vdash_\IL$. Thus:

\ex.[(i)]
$a=b\;$ iff  $\:a\leq^\A b$ and $b\leq^\A a$.

To show that $\A$ is a Heyting algebra, we check that each of the equations defining (the variety of) Heyting algebras is valid in $\A$. This is completely straightforward. Indeed, we have:

\ex.[(ii)]
For each defining equation $s=t$ of Heyting algebras there is a valuation $v$ and formulas $\f,\p$ such that $\overline{v}(\f)=s$, $\overline{v}(\p)=t$, and $\f\vdash_\IL\p$ and $\p\vdash_\IL\f$.\footnote{\label{fn:heyting}For the record, $\A$ is a Heyting algebra iff $(A,0^\A,1^\A,\land^\A,\lor^\A)$ is a bounded lattice and in addition the following equations are valid (omitting superscripts): \vspace{-1ex}
\be
\item[] $a \to a = 1$\vspace{-1ex}
\item[] $a\land (a \to b) = a \land b$\vspace{-1ex}
\item[] $(a\to b) \land b = b$\vspace{-1ex}
\item[] $a\to (b \land c) = (a \to b) \land (a \to c)$\vspace{-1ex}
\ee
To be more accurate, we should have distinguished in (ii) between the syntactic equation $s=t$ in which $s,t$ are \emph{terms}, and the corresponding equality in $\A$, but we trust our abuse of notation here is not a problem.
}

It then follows, from (i) and consistency with $\vdash_\IL$, that $\A$ is a Heyting algebra. As an example, let us check the equation

\ex.[(iii)]
$a \land^\A(a\to^\A b) = a\land^\A b$

Let $v(p)=a$, $v(q)=b$, $\f = p\land(p\to q)$, and $\p = p\land q$. Since we have that $p\land(p\to q) \vdash_\IL p\land q$, $p\land q\vdash_\IL p\land(p\to q)$, and $\overline{v}$ is a homomorphism, the claim follows. As the example shows, we are merely translating familiar laws of $\IL$ into (in)equalities in $\A$, which are valid by the consistency requirement.

If $\A$ is instead consistent with $\vdash_\CL$, the defining equations of Heyting algebras are still valid. It only remains to show that (dropping superscripts) $a \land a' = 0$ and $a \lor a' = 1$. The first is a consequence of (iii) and the fact that $a\land 0=0$. The second identity is the Law of Excluded Middle.
\QED

We may take left-to-right implication of this theorem to say that $\IL$ and $\CL$ are Carnap categorical in algebraic semantics, with respect to the respective classes of algebras. The fact that the implication goes both ways could be seen as a further justification for taking Heyting (Boolean) algebras to be the `standard interpretations'; see also Section \ref{ss:meaning} below.

\begin{remark}
Just as for possible worlds semantics, to enforce that all algebraic interpretations consistent with $\vdash_\IL$ are standard (Heyting), it is not enough to require that all $\IL$-\emph{theorems} are valid. The proof of Theorem \ref{mainthmalg} doesn't work if we cannot conclude from $\p\vdash_\IL\f$ that (for all $v$) $\overline{v}(\p) \leq^\A \overline{v}(\f)$. In fact, the following counter-example can be found in BH19.

If $H$ is a Heyting algebra and $j$ a nucleus on $H$, BH19 defines the algebra $D(H,j)$ of the same signature to be just like the Heyting algebra of fixpoints $H_j$ (defined immediately after Definition \ref{def:nuclear}), except that $a \to^j b =_\text{def} a \to jb$. Algebras of this form are called \emph{Dummett algebras} in BH19.\footnote{Dummett algebras have an interesting connection to a proposal in \cite{dummett00} to explain the significance of Beth semantics in terms of a distinction between a formula being \emph{verified} and \emph{assertible} at a stage $x$ (BH19, end of Sect.\ 3.2).
}
They need not be Heyting algebras --- Bezhanishvili and Holliday give a simple example --- but they of course qualify as algebraic interpretations. However, BH19 also proves that the formulas valid in all Dummett algebras are \emph{exactly} the theorems of $\IL$ (their Theorem 3.25). Thus, by Theorem \ref{mainthmalg}, non-Heyting Dummett algebras need not be consistent with $\vdash_\IL$.
\end{remark}

\section{Discussion\label{s:discussion}}

We find it rather remarkable that $\IL$ is Carnap categorical. Many scholars think that the intuitionistic meaning of the connectives is best captured by the informal Brouwer-Heyting-Kolmogorov explanation, and that the familiar formal semantics for which $\IL$ is sound and complete fail in various degrees to do justice to that explanation. For example, the BHK explanation of the meaning of $\to$ uses notions like `proof' and `construction': a proof of $\f\to\p$ is a construction that takes any proof of $\f$ to a proof of $\p$. Nuclear semantics involves no similar notions. Nor does it rely on facts about verification or assertion, although one can argue quite convincingly, as BH19 does building on Dummett and others, that some instances of nuclear semantics \emph{represent} such facts rather accurately.

\subsection{Carnap's question and the meaning of the connectives\label{ss:meaning}}

However, Carnap's question, as we construe it here, is a model-theoretic question. If you will, it is about the \emph{extension} of the meaning of the connectives. And then, relative to certain set-theoretic objects taken to be \emph{semantic values} of sentences (the `extensional part' of sentence meanings), there is nothing more to say about those extensions than how they determine the semantic value of a complex sentence from the semantic values of its (immediate) constituents. In possible worlds semantics for intuitionistic logic, the values are certain upward closed subsets of some poset, so the connectives must be interpreted as functions on those values. And what we find remarkable is that in every nuclear semantics (and similarly for topological semantics), however the nucleus $j$ selects appropriate semantic values, the consequence relation $\vdash_\IL$ \emph{uniquely fixes}, among all the in principle available options, the functions interpreting the connectives on those values.

Similarly, though much less surprisingly, for algebraic semantics. Facing the extensive literature on algebraic semantics for intuitionistic logic, someone might naively ask: Why Heyting algebras? Aren't there other possible algebraic interpretations? Theorem \ref{mainthmalg} gives the answer: if you want to work with a class $C$ of algebras --- algebraic interpretations in our sense --- all of whose members are consistent with $\vdash_\IL$, then $C \sub \mathit{HA}$. Indeed, $\mathit{HA}$ is the \emph{largest} such class. Similarly for $\vdash_\CL$ and $\mathit{BA}$.

Note that these facts have little to do with completeness. To see this, define for each algebraic interpretation $\A$ the semantic consequence relation $\models_\A$ by

\ex.
$\p_1,\ldots,\p_n\models_\A\f\:$ iff for all valuations $v$ on $A$ and all $c\in A$, $c \leq^\A \overline{v}(\p_i)$ for $i=1,\ldots,n$ implies $c \leq^\A \overline{v}(\f)$.

Thus, $\A$ is \emph{consistent} with $\vdash$ iff $\:\vdash\; \sub\; \models_\A$. And $\vdash$ is \emph{sound and complete} for a class $C$ of algebras iff $\:\vdash \;\sub \bigcap_{\A\in C}\!\models_\A$ (soundness) and $\:\bigcap_{\A\in C} \!\models_\A \;\sub\; \vdash$ (completeness). So while it happens to be true that $\mathit{HA}$ is the largest class for which $\vdash_\IL$ is sound and complete, it is also the largest class for which $\vdash_\IL$ is sound, by Theorem \ref{mainthmalg}. Completeness is not needed to single out $\mathit{HA}$.

We note further that the abstract completeness of a logic, in the sense of its set of theorems (or its consequence relation) being recursively enumerable, is irrelevant to categoricity. For example, every \emph{intermediate logic} (logic between $\IL$ and $\CL$) is categorical with respect to possible worlds semantics as per Theorem \ref{mainthm}, but there are uncountably many intermediate logics, hence uncountably many whose set of theorems is not recursively enumerable.

Going beyond propositional logic, a concrete example of an incomplete logic which is Carnap categorical (the meaning of the logical constants is fixed by the standard consequence relation) is the logic $\mathcal{L}(\Q_0)$, which is classical first-order logic with the additional quantifier `there are infinitely many'.\footnote{This was observed by the second author, and is stated and generalized in \cite{speitel20}. The set of valid sentences in $\mathcal{L}(\Q_0)$ is not recursively enumerable.}


\subsection{Other logics, other settings}

Here are some final observations and open questions.\footnote{Questions and comments from several people inspired the remarks in this subsection, in particular from Wes Holliday for parts 1 and 2, from Johan van Benthem for part 3, and from Denis Bonnay for part 4.
}
\vspace{1ex}

\noindent
\textbf{1.} An obvious question is if our results here are best possible for intuitionistic propositional logic. One way of making this precise is to let, for $\Phi\sub\{\neg,\land,\lor,\to\}$, $L_\Phi$ be the propositional language with connectives in $\Phi$, and $\vdash_\Phi$ be the consequence relation defined by the usual natural deduction rules for these connectives. Here we use $\neg$ rather than $\bot$,\footnote{Say, with the rules
\bi
\item[] 
\alwaysNoLine
\AxiomC{$\p$ $\neg\p$} 
\alwaysSingleLine
\UnaryInfC{$\f$}
\DisplayProof
\hspace{10ex}
\alwaysNoLine
\AxiomC{$[\f]$} 
\UnaryInfC{$:$} 
\UnaryInfC{$\p$} 
\AxiomC{$[\f]$} 
\UnaryInfC{$:$} 
\UnaryInfC{$\neg\p$} 
\alwaysSingleLine
\BinaryInfC{$\neg \f$} 
\DisplayProof
\ei
}
so, modulo this change (and deleting curly brackets and commas), $\vdash_{\neg\land\lor\to}\;=\;\vdash_\IL$. 
Now we can ask (with the obvious modification of Definition \ref{def:categoricity}): If $\Phi$ is a proper subset of $\{\neg,\land,\lor,\to\}$, is $\vdash_\Phi$ categorical with respect to some suitable semantics?

Again we must distinguish algebraic semantics from possible worlds semantics. Once the relevant classes of algebras have been identified, algebraic categoricity seems (again) essentially obvious. For example, in analogy to Theorem \ref{mainthmalg} we have:\footnote{The signature of $\A$ in (a) is $\{\land,\lor\}$, and in (b) it is $\{\neg,\land,\lor\}$. So there is no 1, but as per footnote \ref{fn:consistency}, consistency is still well-defined. The distributive lattices in (a) need not be bounded, for example, any linearly ordered set is consistent with $\vdash_{\land\lor}$; note that this logic has no theorems. The lattices in (b) are bounded, and pseudo-complementation means that for each $a\in A$, the maximum of $\{b\!: a\land b=0\}$ belongs to $A$. See \cite{fontverdu91} for $\vdash_{\land\lor}$, and  \cite{rebagliatoverdu93} for $\vdash_{\neg\land\lor}$, or consult \cite{font16}. It is well-known that $\vdash_{\land\lor}$ can also be defined as $\:\vdash_\IL \restriction \!L_{\land\lor}$, and similarly for $\vdash_{\neg\land\lor}$.
}

\begin{fact}
\bi
\item[\rm (a)] $\A$ is consistent with $\vdash_{\land\lor}$ iff $\A$ is a distributive lattice.\vspace{-1ex}
\item[\rm (b)] $\A$ is consistent with $\vdash_{\neg\land\lor}$ iff $\A$ is a pseudo-complemented distributive lattice.
\ei
\end{fact}

Possible worlds categoricity is a different matter. The case of $L_\emp$ is trivial (we have $\Gamma\!\vdash_\emp\!\f \Leftrightarrow \f\in\Gamma$, but there are no connectives to interpret), and $\vdash_\land$ is easily seen to be categorical ($I^\F(\land)$ must be intersection), but what about others? A particularly relevant case is $\vdash_{\neg\land\lor}$, since the proof we gave that $\lor$ must be standard (in every intuitionistic interpretation of the relevant kind consistent with $\vdash_\IL$) depended on facts about $\to$. One may conjecture that this is necessary, and so that $\vdash_{\neg\land\lor}$ is \emph{not} Carnap categorical with respect to, say, Kripke semantics, but this is open. We can, however, give an example where there is algebraic categoricity but not possible worlds categoricity. Consider the logics $\vdash_\lor$ and $\vdash_\neg$. First, it is easy to see that

\ex.\label{algor}
$\A$ is consistent with $\vdash_{\lor}$ iff $\A$ is a semilattice.\footnote{Now the signature is $\{\lor\}$, and we define $a\leq^\A b \Leftrightarrow a \lor^\A b = b$; see footnote \ref{fn:consistency}. We refrain from fomulating a corresponding fact for $\vdash_\neg$; the algebras of that signature are not algebraic interpretations in our sense.
}

Second, however:

\begin{fact}
Neither $\vdash_\lor$ nor $\vdash_\neg$ is Carnap categorical with respect to Kripke semantics.
\end{fact}

\proof (outline) 
$\vdash_\lor$ proves that disjunction is commutative, associative, and idempotent. Using this, it is not hard to verify: 

\ex.\label{a1}
$\p_1,\ldots,\p_n \vdash_\lor \f$ iff there is $\p_i$ such that each atom in $\p_i$ occurs in $\f$.

Now let $\F$ be any poset such that there is a nucleus $j$ on $\mathit{Up}(\F)$ which is not the identity, and let $I^\F(U,V)=j(U\cup V)$. Using \ref{a1} and the monotonicity of $j$, it is easy to see that $I^\F$ is consistent with $\vdash_\lor$, but $I^\F$ is not the standard Kripke interpretation.\footnote{This simple non-standard interpretation was suggested by Wesley Holliday, and replaces a more \emph{ad hoc} construction that we originally gave.
}

Next, consider $\vdash_\neg$. As is well-known:

\ex.\label{neg1}
\a. $\f,\neg\f\vdash_\neg\p$
\b. $\f\vdash_\neg \neg\neg\f $
\b. $\neg\f \dashv\vdash_\neg \neg\neg\neg\f$

Every $L_\neg$ formula is of the form $\neg^n p$ for some $n\geq0$, where $\neg^0p=p$ and $\neg^{n+1}p= \neg\neg^n p$. Say that $\{\neg^{n_1}p_1,\ldots,\neg^{n_k}p_k\}$ is \emph{contradictory} if $p_i = p_j$ for some $i,j$, and $|n_i-n_j|$ is odd. Then one can show:

\ex.\label{neg2}
$\neg^{n_1}p_1,\ldots,\neg^{n_k}p_k \vdash_\neg \neg^{m}q\;$ iff either $\{\neg^{n_1}p_1,\ldots,\neg^{n_k}p_k\}$ is contradictory, or ($q = p_i$ for some $i$ and ($n_i,m$ are both odd, or $n_i,m$ are both even and $m\geq2$)).

Using this, one can verify that the rules (R1) $\f\vdash\neg\neg\f$, (R2) $\neg\neg\neg\f\vdash\neg\f$, and (R3) $\f,\neg\f\vdash\p$, provide a Hilbert style axiomatization of $\vdash_\neg$.

Now let $\F$ be the 4-element poset $(X,\leq)$, where $X=\{0,0',a,b\}$, $0<a$, and $0'<b$. The upsets are $\emp,X,\{a\},\{0,a\},\{b\},\{0',b\}$, and $\{a,b\}$. Recall that $I^\F_\text{st}(\neg)(U) = \{x\!: \forall y\geq x\: y\not\in U\}$, and define the interpretation $I^\F(\neg) = N$ as follows:

\ex.\label{neg3}
$N(U)=I^\F_\text{st}(\neg)(U)$, \emph{except} that $N(\{a\})=\{b\}$ and $N(\{b\})=\{a\}$. 

(Note that $I^\F_\text{st}(\neg)(\{a\})=\{0',b\}$.) Then, for all upsets $U$,

\ex.\label{neg4}
\a. $U \sub N(N(U))$
\b. $N(N(N(U))) \sub N(U)$
\c. $U\cap N(U)=\emp$

Using \ref{neg4}, an induction over the length of Hilbert style proofs shows that $I^\F$ is consistent with $\vdash_\neg$.
\QED

Thus, in contrast with the case of $\land$, the introduction and elimination rules for $\lor$ do \emph{not} fix its meaning in Kripke semantics. This is also in contrast with classical 2-valued semantics, where consistency with those rules does fix the meaning, i.e.\ the standard truth table, for $\lor$. Similarly for $\neg$.\footnote{The introduction rules for $\lor$ fix the first three rows, (1,1,1), (1,0,1), and (0,1,1), in the truth table for $\lor$. The fourth row, (0,0,0), is fixed by the fact that $p\lor p\vdash_\lor p$, which is an instance of the elimination rule for $\lor$. Likewise, $p,\neg p\vdash_\neg q$ and $p\vdash_\neg \neg\neg p$ together fix the standard truth table for $\neg$. Note that $\vdash_\lor$ is the restriction of $\vdash_\IL$, as well as of $\vdash_\CL$, to $L_\lor$. This of course fails for $\vdash_\neg$, which nevertheless fixes the 2-valued meaning of $\neg$.
}
\vspace{2ex}

\noindent
\textbf{2.} One can also weaken the logics $\vdash_\Phi$ by constraining the inference rules in various ways. A case in point is \cite{holliday22}, which studies a weaker logic in $L_{\neg\land\lor}$, called $\vdash_\mathsf{F}$, with restrictions on the $\lor$E and $\neg$I rules. Among other things, distributivity no longer holds.\footnote{Holliday presents $\vdash_\mathsf{F}$ with a Fitch style natural deduction system, where the added constraint becomes a requirement on the Reiteration rule. $\vdash_\mathsf{F}$ is also extended to a first-order language with the quantifiers $\forall$ and $\exists$ (but without $\to$ and $=$). Failure of distributivity is a characteristic of quantum logic, but Holliday argues that it also accords with certain facts about natural language semantics, in particular facts about epistemic modals.
}
He shows how $\vdash_\mathsf{F}$ can be seen as a neutral base logic, from which intuitionistic logic, and versions of the \emph{orthologic} studied in \cite{goldblatt74}, and classical logic, can be obtained by suitable additions or changes to the rules, or by corresponding constraints in the algebraic or the relational semantics he provides for $\vdash_\mathsf{F}$. Finding a non-standard interpretation consistent with $\vdash_\mathsf{F}$ --- if there is such an interpretation --- might be easier than finding one for $\vdash_{\neg\land\lor}$.

In the same spirit, we may ask if there is a proper fragment of $\vdash_\IL$ in the language $L_{\neg\land\lor\to}$ which is Carnap categorical. If no such fragment exists, that would be a new kind of functional completeness property of $\IL$.\vspace{2ex}

\noindent
\textbf{3.} Theorem \ref{mainthm} is an existence and uniqueness result: \emph{there is} an interpretation consistent with $\vdash_\IL$ (namely, the standard interpretation), and it is in fact the \emph{only} one. In the proof-theoretic tradition, existence and uniqueness of propositional connectives is a well-established topic; see \cite{humberstone11}, Ch.\ 4, for a comprehensive overview. Here, existence and uniqueness is relative to a set of \emph{rules}. For example, $\land$ and $\lor$ are \emph{unique} relative to a natural deduction presentation of $\IL$, in the sense that if we introduce new connectives $\land'$ and $\lor'$, with the same rules as for $\land$ and $\lor$, respectively, then $\f\land\p$ is equivalent to $\f\land'\p$, and similarly for $\lor$ and $\lor'$. Indeed, we only need the introduction and elimination rules for these two connectives. With the notation from part 1 above, we have: $\f\land\p\dashv\vdash_{\land\land'}\f\land'\p$ and $\f\lor\p\dashv\vdash_{\lor\lor'}\f\lor'\p$.

As Humberstone noted, and \cite{dosenSH88} spells out in detail, there is a difference between $\land$ and $\lor$ here: while the former is \emph{implicitly definable} in a precise sense, the latter, although unique, is not. Do\v{s}en and Schroeder-Heister explore connections to Beth's Definability Theorem.\footnote{\cite{humberstone11}, Chapter 4.35, proves a similar result about $\lor$: it is not uniquely characterizable by what he calls \emph{zero-premiss rules}. See also the discussion in \cite{williamson88}, pp.\ 110--114.
}
It is an intriguing question if, and in that case how, their proof-theoretic notion of implicit definability relates to our semantic notion of categoricity (which, as we have seen, exhibits a similar contrast between $\land$ and $\lor$).\vspace{2ex}

\noindent
\textbf{4.} Can Theorem \ref{mainthm} be generalized to arbitrary intuitionistic interpretations as specified in Definition \ref{def:possworldssem}? This concerns in particular the function $I^{\F,\mathit{at}}$, which could be seen as the interpretation of an invisible logical constant preceding all atoms. Does categoricity extend to that constant? More precisely, if $I^\F$ is consistent with $\vdash_\IL$, under what circumstances would $I^{\F,\mathit{at}}$ have to be a nucleus?
\vspace{2ex}

\noindent
\textbf{5.} We end with another set of questions, not related specifically to intuitionistic logic but the general setting in which Carnap's question has been posed. We have so far simply followed the lead of \cite{bonwes16} and other papers in considering consequence relations $\vdash$ where the first argument is a (possibly empty) \emph{set} of formulas and the second is a \emph{formula}; the \textsc{set-fmla} setting in the terminology of \cite{humberstone11}. One could look at other settings, like \textsc{fmla} (no premises), \textsc{set-set}, \textsc{fmla-fmla}, or \textsc{seq-fmla} (with a sequence of premises). As we have seen, the categoricity results would differ significantly. We think \textsc{set-fmla} most naturally corresponds to intuitive ideas about `what follows from what', but one would like a better argument for why the others are less suitable. Or, less contentiously, an overview of what the results would be in those settings.\footnote{In fact, the solution proposed in \cite{carnap43} to the existence of non-normal interpretations was in effect to use the \textsc{set-set} framework instead. This makes good sense for classical logic. But for intuitionistic logic there is, as far we know, no obviously correct candidate for the \textsc{set-set} version of $\IL$. For example, $p\lor q\vdash p,q$ is not valid in standard Beth semantics (since we can have $j_b(U\cup V) \not\sub U\cup V$), although it holds in Kripke semantics.
}

In this context, here is final more concrete issue. In all categoricity results (in our sense) that we are aware of, it is practically immediate that any interpretation consistent with the relevant consequence relation $\vdash$ must interpret conjunction ($\land$) standardly. It is the other logical constants that need some work (except for $\bot$ if \emph{Ex Falso} holds). But one could argue that this is built into our definition of consistency with $\vdash$ in Definition \ref{def:possworldssem}:4: the \emph{intersection} of the values of the premises must be included in the value of the conclusion. This guarantees that $U\cap V \sub I(\land)(U,V)$; the converse inclusion is just from single-premise $\land$-elimination. 

This choice should at least be motivated. Is it necessary? One approach could be to see what happens in the \textsc{fmla-fmla} setting: does it too force $\land$ to be standard? The answer seems to be No in a general possible worlds environment; at least if interpretations now have the form $\F = (X,\C,\sub)$, where $\C$ is any subset of $\P(X)$, and semantic values are required to stay in $\C$. It is easy to find such $\C$ such that $I^\F$ is `\textsc{fmla-fmla}-consistent' with $\vdash$ --- in the sense that if $\p\vdash\f$ then for all valuations $v$, $\lint\p\rint^{I^\F}_v \sub \lint\f\rint^{I^\F}_v$ --- but $I^\F(\land)$ is \emph{not} intersection. On the other hand, if $\C$ is required to be closed under intersections, then $I^\F(\land)$ will be intersection.\vspace{3ex}

\noindent
We feel these issues deserve more thought, but leave them for future work.

\subsection*{Acknowledgements}
We thank Valentin Goranko and Lloyd Humberstone for helpful comments at a (much) earlier stage of this work, and Johan van Benthem, Denis Bonnay, Valentin Goranko (again), and Wesley Holliday for constructive remarks on the present paper. Special thanks to Wes for several helpful conversations and for allowing us to present his Example \ref{ex:holliday}.

\printbibliography[heading=bibintoc]

\end{document}